\documentclass[11pt, twoside]{article}
\usepackage{fullpage,ifthen,epsf,hyperref,url,breakurl,enumerate,comment,flushend,multirow}
\usepackage[active]{srcltx}
\usepackage{amsthm,amsfonts,amsmath,amssymb}
\usepackage{microtype,multicol,sidecap,mathrsfs}
\usepackage{esint}
\usepackage{algorithm}
\usepackage[noend]{algpseudocode}
\usepackage[lofdepth,lotdepth]{subfig}
\usepackage[usenames,dvipsnames,svgnames,table]{xcolor}
\usepackage{graphicx}
\graphicspath{ {./figures1/} {./figures2/} }
\usepackage[normalem]{ulem}

\useunder{\uline}{\ul}{}

\newlength{\widebarargwidth}
\newlength{\widebarargheight}
\newlength{\widebarargdepth}

\makeatletter
\long\def\@makecaption#1#2{
        \vskip 0.8ex
        \setbox\@tempboxa\hbox{\small {\bf #1:} #2}
        \parindent 1.5em 
        \dimen0=\hsize
        \advance\dimen0 by -3em
        \ifdim \wd\@tempboxa >\dimen0
                \hbox to \hsize{
                        \parindent 0em
                        \hfil 
                        \parbox{\dimen0}{\def\baselinestretch{0.96}\small
                                {\bf #1.} #2
                                } 
                        \hfil}
        \else \hbox to \hsize{\hfil \box\@tempboxa \hfil}
        \fi
        }
\makeatother

\begin{document}
\title{Impact of Detour-Aware Policies on Maximizing Profit in Ridesharing}

\author{
Arpita Biswas$^{1,\dagger}$, Ragavendran Gopalakrishnan$^{2}$, Theja Tulabandhula$^{3}$,\\
Asmita Metrewar$^{4}$, Koyel Mukherjee$^{5}$ and Raja Subramaniam Thangaraj$^{2}$\\
\begin{tabular}{c}
$^{1}$Indian Institute of Science\\
$^{2}$Conduent Labs India\\
$^{3}$University of Illinois Chicago\\
$^{4}$Apple India\\
$^{5}$IBM Research India\\
$^{\dagger}$Email: arpita.biswas@csa.iisc.ernet.in
\end{tabular}
}
\date{May 4, 2017}

\maketitle
\begin{abstract}
This paper provides efficient solutions to maximize profit for commercial ridesharing services, under a pricing model with detour-based discounts for passengers. We propose greedy heuristics for real-time ride matching that offer different trade-offs between optimality and speed. Simulations on New York City (NYC) taxi trip data show that our heuristics are up to $90\%$ optimal and $10^5$ times faster than the (necessarily) exponential-time optimal algorithm.

Commercial ridesharing service providers generate significant savings by matching multiple ride requests using heuristic methods. The resulting savings are typically shared between the service provider (in the form of increased profit) and the ridesharing passengers (in the form of discounts). It is not clear a priori how this split should be effected, since higher discounts would encourage more ridesharing, thereby increasing total savings, but the fraction of savings taken as profit is reduced. We simulate a scenario where the decisions of the passengers to opt for ridesharing depend on the discount offered by
the service provider. We provide an adaptive learning algorithm \texttt{IDFLA} that learns the optimal profit-maximizing discount factor for the provider. An evaluation over NYC data shows that \texttt{IDFLA}, on average, learns the optimal discount factor in under $16$ iterations.

Finally, we investigate the impact of imposing a detour-aware routing policy based on sequential individual rationality, a recently proposed concept. Such restricted policies offer a better ride experience, increasing the provider's market share, but at the cost of decreased average per-ride profit due to the reduced number of matched rides. We construct a model that captures these opposing effects, wherein simulations based on NYC data show that a $7\%$ increase in market share would suffice to offset the decreased average per-ride profit.
\end{abstract}

\section{Introduction}

Ridesharing is a key initiative that can curb ever-increasing congestion on the urban transportation network and is an important means to move the world towards a sustainable future. While the term \textit{ridesharing} has been used to include peer-to-peer carpooling platforms as well, here we focus on commercial ridesharing service providers such as LyftLine and UberPool, that hold a significant share of the ridesharing population.

While ridesharing is undoubtedly appealing from a sustainability perspective, profit maximization is the primary goal of commercial providers. Given the emerging literature on detour-aware pricing and routing policies, e.g.,~\cite{Jung2013,DBLP:journals/corr/GopalakrishnanM16} that enhance the ride experience, algorithms for profitable ride-matching should adapt to keep up. Such quality-enhancing policies encourage increased adoption of ridesharing, but perhaps at the cost of reduced average per-ride profit. We study three connected issues centered around profit optimization in commercial ridesharing.

First, it is unclear how to tune real-time ride-matching algorithms to quickly find a near-optimal set of matches that results in high profit under detour-based discounts (the exact optimization is an NP-hard problem~\cite{Cordeau2007}).

Second, ridesharing generates significant economic savings by reducing the driver-miles necessary to serve a set of passengers, which are split between the service provider (their profit) and the passengers (as discounts). It is not obvious how this split should be effected, since higher discounts encourage more ridesharing, increasing total savings, but the fraction of savings that constitutes the profit is smaller.

Third, in order to survive in an increasingly competitive market, ridesharing providers adopt additional quality-enhancing features which do not align fully with the goal of profit maximization, perhaps resulting in lower profit (due to the smaller feasible set of quality-compliant rides). It is important to understand the overall impact on profitability, and ensure that any resulting increase in market share can counter the reduced average per-ride profit.

The issues mentioned above are crucial for commercial ridesharing platforms; however, they are less explored. We address these issues, give solutions, and provide extensive evaluation on a real-world dataset.

\subsection{Our Contributions}
We investigate three important aspects of profit optimization in commercial ridesharing, resulting in the following contributions:
\begin{itemize}
\item In Section~\ref{sec:8}, given a detour-based discount policy, we provide an Integer Linear Program (ILP), and a family of efficient \textit{greedy} heuristics that quickly match passengers to rides to approximately maximize profit. We evaluate their performance experimentally.
\item In Section~\ref{sec:learning}, under a model in which the likelihood of opting for ridesharing depends on the offered discount, we propose an adaptive algorithm that learns the profit-maximizing discount parameter.
\item In Section~\ref{sec:qos-impact}, we additionally impose a detour-aware routing policy based on \textit{sequential individual rationality}, introduced in~\cite{DBLP:journals/corr/GopalakrishnanM16}. Using simulations, we investigate the minimum increase in market share needed to offset the reduced average per-ride profit.
\end{itemize}

\section{Related Work}
There is a huge body of literature that studies several optimization problems related to ride sharing~\cite{Furuhata2013,Agatz2012,Pelzer2015}, but the problem of pricing in ridesharing has been relatively less studied, with profit optimization rarely discussed. \cite{Banerjee2015} characterizes the pros/cons of static versus dynamic pricing. \cite{Jung2013} provides heuristics for maximizing profit subject to bounded detours. In contrast, we study profit maximization over a pricing model that incorporates a detour-based discount.

\cite{nguyen} proposes a dynamic, demand-responsive ridesharing system using auction mechanisms, where customers enjoy additive detour-based discounts. Our model considers \textit{multiplicative} detour-based discounts; moreover, we also provide an algorithm to learn the profit-maximizing discount parameter under a suitable user choice model.

A closely related problem is that of cost sharing in peer-to-peer ridesharing or carpooling, which, until recently, received little attention. Individual passengers are either asked to post what they are willing to pay in advance~\cite{sharek}, to share the total cost proportionately among themselves according to the distances travelled~\cite{Geisberger2010,Agatz2011}, or negotiate their shares on their own during/after the ride. Such methods ignore the real-time costs and delays incurred during the ride and are generally insensitive to the disproportionate delays encountered during the ride. Fair cost sharing in ridesharing has been studied~\cite{bistaffa,Kleiner2011,zhao,kamar,DBLP:journals/corr/GopalakrishnanM16}, but its impact on profitability is left unexplored. We fill this gap by studying the trade-off between average per-ride profits and market share.

\section{Matching Rides to Maximize Profit}\label{sec:8}

Ridesharing generates significant profit to commercial platforms such as Uber and Lyft. The profit from a single cab is the difference between the fare paid by the passengers and the amount paid to the driver. When compared to assigning each user their own cab, ridesharing requires fewer cabs to serve the same set of users, leading to increased profit.

At the core of these platforms is a real-time matching algorithm that runs at regular time intervals, and generates a set of matched users from a set of waiting users (preprocessed to accommodate other spatio-temporal constraints). Such a system would maintain a pool of pending requests, and then invoke the matching algorithm every few minutes to determine optimal matches. We provide matching heuristics that directly optimize profit, under a detour-based discount scheme. Any unmatched requests left at the end of the algorithm either wait for another round, or, after a certain waiting time threshold elapses, are assigned fresh, empty cabs. Thus, it is important to provide a time-efficient matching algorithm for the larger dynamic ridesharing system to be sufficiently responsive to incoming user requests.

Moreover, under a detour-based discount policy, shared rides with smaller detours are more likely to be profitable. However, such a myopic observation could lead to suboptimal solutions, thus rendering the problem nontrivial.

In Section~\ref{sec:ILP}, we provide an Integer Linear Program (ILP) that matches $n$ users to shared cabs (each with capacity $\zeta$) to maximize profit, given a multiplicative detour-based discount $f_p(\delta_i,\tau_i)$, which is a \textit{linear} function of the distance-wise ($\delta_i$) and time-wise ($\tau_i$) fractional detours experienced by user $i$. $c_b$ is the ``base'' (fixed) cost for a ride, $c_d$ and $c_t$ are the costs per unit distance and time, respectively. $f_d$ is the fraction of driver earnings taken by the service provider.

\subsection{ILP Formulation}\label{sec:ILP}

Each user $i$ is associated with a source-destination pair $(S_i,D_i)$. We consider $n$ initially empty cabs and construct a complete directed graph $G=(V,E)$ with $2n$ vertices that represent all the sources and destinations, that is, $V=\{S_i,D_i\}_{i=1}^n$. Each directed edge $e=(u,v)$ represents the best route from location $u$ to $v$, and weights $w_d(e)$ and $w_t(e)$ denote the corresponding distance and time, respectively.

The ILP seeks values to assign to a set of optimization variables in order to maximize an objective function subject to a set of constraints.\\

\noindent\textbf{Optimization Variables:} Let $x_{i,j,e}\in\{0,1\}$ denote whether or not user $i$ is served by cab $j$ along edge $e$. Let $y_{j,e}\in\{0,1\}$ denote whether cab $j$ travels along edge $e$ while serving one or more users. Let $z_j\in\{0,1\}$ denote whether cab $j$ serves at least one user. We set $x_{i,j,(u,v)}=0$ if $(u,v)\notin E$.\\

\noindent\textbf{Objective Function:}\\
For each user $i$, we have:
\begin{enumerate}[(a)]
\item Distance travelled, $d_i = \sum_{j=1}^n\sum_{e\in E}x_{i,j,e}w_d(e)$.
\item Time in travel, $t_i = \sum_{j=1}^n\sum_{e\in E}x_{i,j,e}w_t(e)$.
\item Distance-wise fractional detour, $\delta_i = \frac{d_i-w_d(S_i,D_i)}{w_d(S_i,D_i)}$.
\item Time-wise fractional detour, $\tau_i = \frac{t_i-w_t(S_i,D_i)}{w_t(S_i,D_i)}$.
\item The passenger's fare is computed as,
 $$C_i\!=\!(1\!\!-\!\!f_p(\delta_i,\!\tau_i))(c_b + c_dw_d(S_i,\!D_i)\! +\! c_tw_t(S_i,\!D_i)).$$
\end{enumerate}

\noindent For each cab $j$, we have:
\begin{enumerate}[(a)]
\item Distance traveled, $\sigma_j = \sum_{e\in E} y_{j,e}w_d(e)$.
\item Time in travel, $\eta_j = \sum_{e\in E} y_{j,e}w_t(e)$.
\item Driver earnings, $E_j = (1-f_d)(c_bz_j + c_d\sigma_j + c_t\eta_j)$.
\end{enumerate}

\noindent Thus, the objective function to be maximized is the total profit, given by
$p = \sum_{i=1}^nC_i - \sum_{j=1}^nE_j$.\\

\noindent\textbf{Constraints:}
\begin{enumerate}[(i)]
\item $y_{j,e}=1$ if and only if $x_{i,j,e}=1$ for some $i$:
    \begin{equation*}
    \begin{split}
    y_{j,e} &\geq x_{i,j,e} \quad \forall i,j\in\{1,2,\ldots,n\}\quad \forall e\in E\\
    y_{j,e} &\leq \sum_{i=1}^n x_{i,j,e} \quad \forall j\in\{1,2,\ldots,n\} \quad \forall e\in E
    \end{split}
    \end{equation*}
\item $z_j=1$ if and only if $y_{j,e}=1$ for some $e$:
    \begin{equation*}
    \begin{split}
    z_j &\geq y_{j,e} \quad \forall j\in\{1,2,\ldots,n\} \quad \forall e\in E\\
    z_j &\leq \sum_{e\in E} y_{j,e} \quad \forall j\in\{1,2,\ldots,n\}
    \end{split}
    \end{equation*}
\item Each user is picked up by exactly one cab:
    \begin{equation*}
    \sum_{u\in V}\sum_{j=1}^n x_{i,j,(S_i,u)} = 1 \quad \forall i\in\{1,2,\ldots,n\}
    \end{equation*}
\item Each user is dropped off by exactly one cab:
    \begin{equation*}
    \sum_{u\in V}\sum_{j=1}^n x_{i,j,(u,D_i)} = 1 \quad \forall i\in\{1,\ldots,n\}
    \end{equation*}
\item Each user is not served before pickup:
    \begin{equation*}
    \sum_{u\in V}\sum_{j=1}^n x_{i,j,(u,S_i)} = 0 \quad \forall i\in\{1,\ldots,n\}
    \end{equation*}
\item Each user is not served after dropoff:
    \begin{equation*}
    \sum_{u\in V}\sum_{j=1}^n x_{i,j,(D_i,u)} = 0 \quad \forall i\in\{1,\ldots,n\}
    \end{equation*}
\item Each user is served by a single cab, only between pickup and dropoff:
    \begin{equation*}
    \sum_{u\in V} x_{i,j,(u,v)} = \sum_{u\in V} x_{i,j,(v,u)} \quad \forall v\in V\setminus\{S_i,D_i\} \ \forall i\ \forall j
    \end{equation*}

\item Each cab serves at most $\zeta$ users along any edge:
    \begin{equation*}
    \sum_{i=1}^n x_{i,j,e} \leq \zeta \quad \forall j\in\{1,\ldots,n\} \quad \forall e\in E
    \end{equation*}
\item Each cab's overall route is not disjoint, i.e., there are no gaps between serving users during which the cab is empty. Equivalently, if a cab serves a total of $k\geq 1$ users, then it must do so using exactly $2k-1$ edges:
    \begin{equation*}
    \sum_{e\in E}y_{j,e} \geq 2\left(\sum_{i=1}^n\sum_{u\in V}x_{i,j,(S_i,u)}\right)-1
    \end{equation*}
    The inequality ($\geq$) is to allow empty cabs ($k=0$).
\item Binary integer constraints:
    \begin{equation*}
    \begin{split}
    x_{i,j,e} &\in \{0,1\} \quad \forall i,j\in\{1,\ldots,n\} \quad \forall e\in E\\
    y_{j,e} &\in \{0,1\} \quad \forall j\in\{1,\ldots,n\} \quad \forall e\in E\\
    z_j &\in \{0,1\} \quad \forall j\in\{1,\ldots,n\}
    \end{split}
    \end{equation*}
\end{enumerate}

Computing the optimal solution to the ILP takes a lot of time and memory, and is unsuitable for real-time matching, which requires fast, near-optimal techniques. However, it is a handy benchmark against which to experimentally evaluate the performance of our proposed methods (Section~\ref{sec:profit-method-comparison}).

\subsection{Heuristics for Real-Time Ride Matching}\label{sec:model}

Given a set $\mathcal{S}$ of users matched to a cab, the driving distance and time for the best route for serving all the users $i\in \mathcal{S}$ are denoted by $d(\mathcal{S})$ and $t(\mathcal{S})$, respectively. For each user $i\in \mathcal{S}$, $d_i(\mathcal{S})$ and $t_i(\mathcal{S})$ denote the driving distance and time from $S_i$ to $D_i$ along that route, respectively. Thus, the profit from a cab serving a set of users $\mathcal{S}$ can be computed as follows:
\begin{itemize}
\item \textit{User Cost:} The cost of the ride to user $i\in\mathcal{S}$,
\begin{equation}
    C_i(\mathcal{S})\!\! = \!\!\left(1\!\!-\!\!f_p(\delta_i(\mathcal{S}),\tau_i(\mathcal{S}))\right)\left(c_b\! +\! c_d d_i(\{i\}) \!+\! c_t t_i(\{i\})\right),\label{eq:passenger_Cost}
\end{equation}

    where $\delta_i(\scriptsize{\mathcal{S}})$=$\frac{d_i(\mathcal{S}) - d_i(\{i\})}{d_i(\{i\})}$, and $\tau_i(\scriptsize{\mathcal{S}})$=$\frac{t_i(\mathcal{S}) - t_i(\{i\})}{t_i(\{i\})}$.
\item \textit{Driver Earnings:} The earnings to the driver from the ride,
\begin{equation}
    E(\mathcal{S}) = \left(1-f_d\right)\cdot\left(c_b + c_d d(\mathcal{S}) + c_t t(\mathcal{S})\right).
\end{equation}

\item \textit{Profit:} The profit to the service provider,
\begin{equation}
    p(\mathcal{S}) = \displaystyle\sum_{i\in\mathcal{S}}C_i(\mathcal{S}) - E(\mathcal{S}).
\end{equation}

\end{itemize}
The incremental profit of \textit{merging} two cabs (combining the users matched to two cabs into one and removing the other cab from the system), serving sets of users $\mathcal{S}_j$ and $\mathcal{S}_k$ is:

\begin{equation}
\Delta p(\mathcal{S}_j,\mathcal{S}_k) = p(\mathcal{S}_j\cup\mathcal{S}_k) - p(\mathcal{S}_j) - p(\mathcal{S}_k).
\end{equation}

The above framework for real-time ride matching is more general than that of the ILP. A limitation of the ILP is that $f_p(\cdot)$ is restricted to be a linear function. Moreover, it may not be possible for an ILP to explicitly incorporate additional detour-aware routing constraints and user preferences. However, during performance evaluation of our proposed methods, we limit the scope to the constraints listed in Section~\ref{sec:ILP}, in order to compare against the optimal ILP solution.\\

\begin{figure*}[!htb]
    \centering
    \subfloat[width=0.45\textwidth][Total profit obtained by various methods \label{fig:profit_20}]{\includegraphics[width=0.45\textwidth]{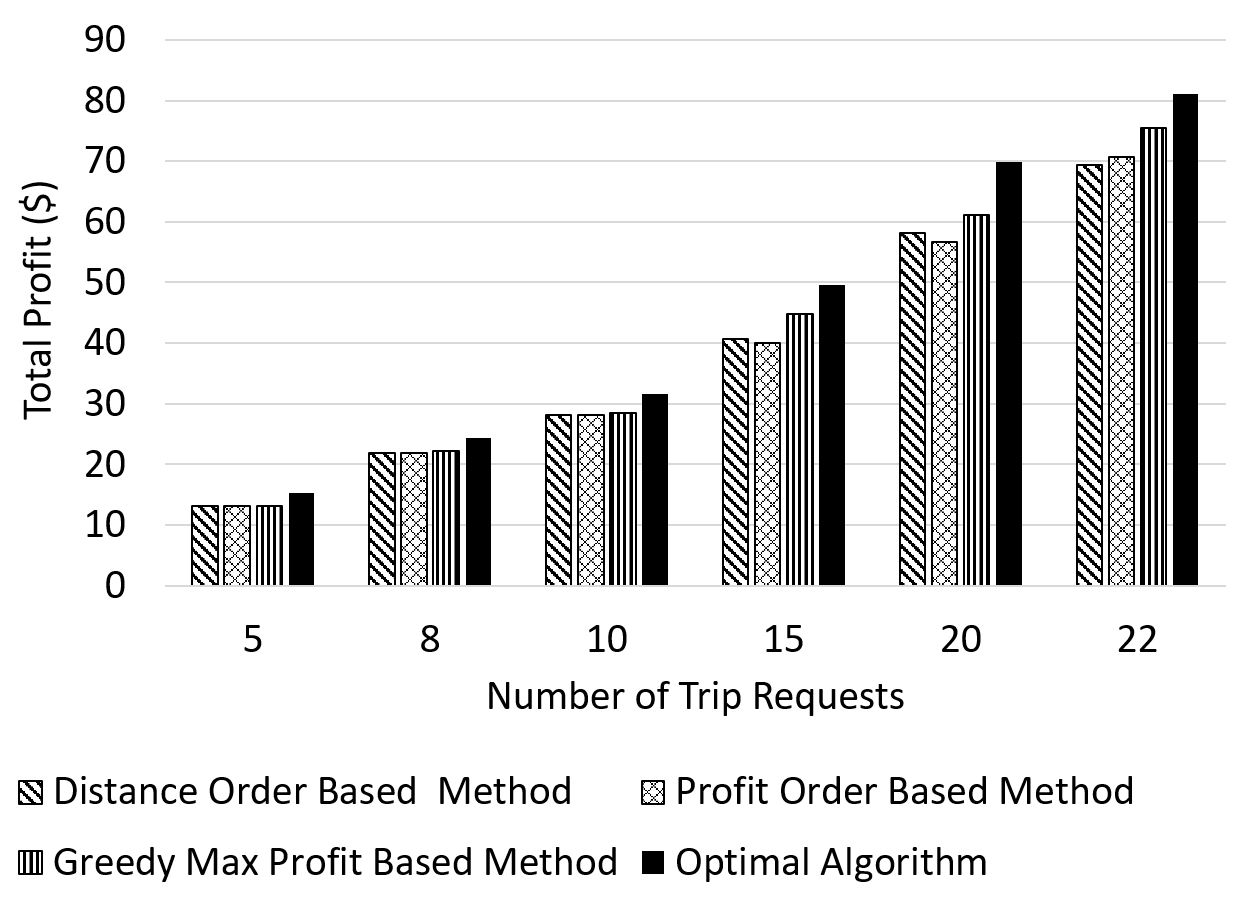}}
    \hfill
    \subfloat[width=0.45\textwidth][Time taken by different methods.\label{fig:time_20}]{\includegraphics[width=0.45\textwidth]{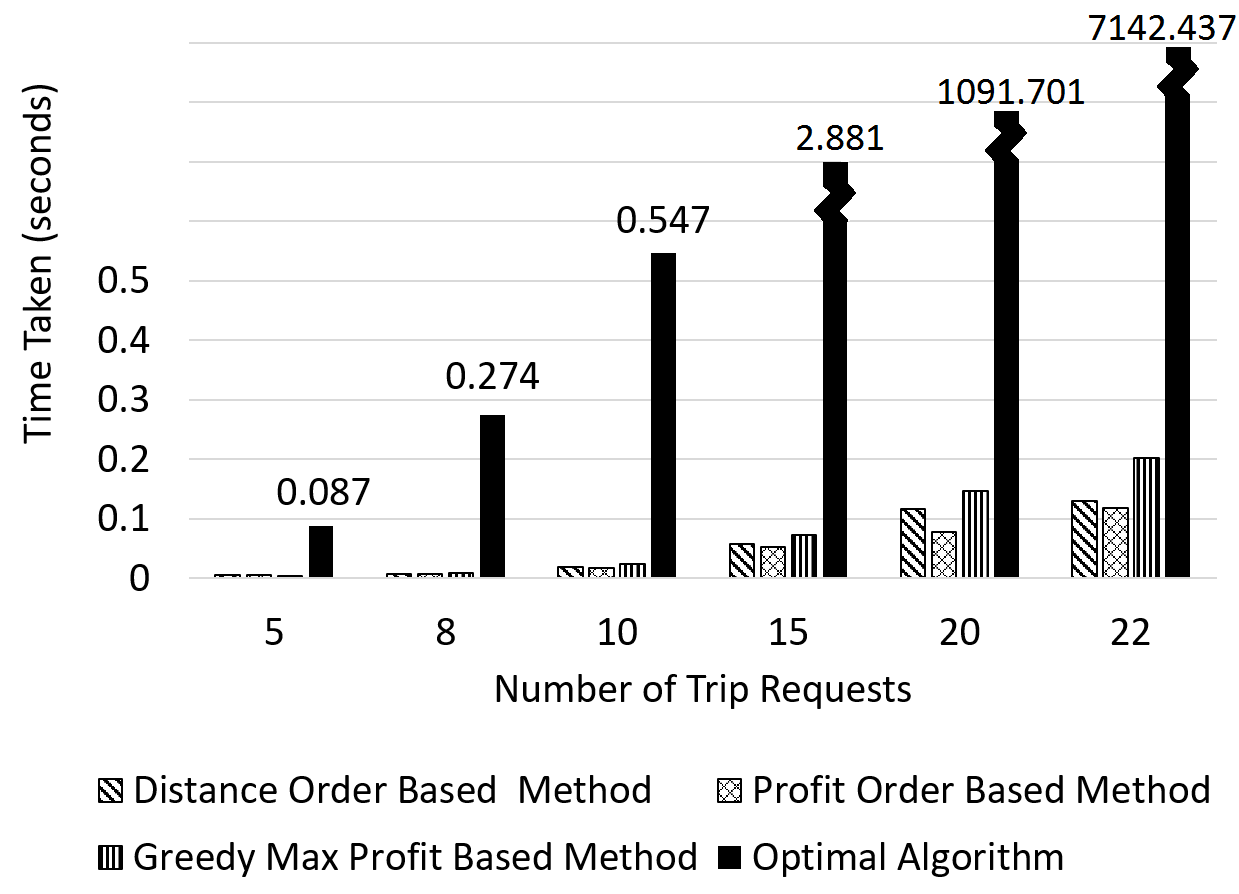}}
\caption{Comparison of various heuristic matching methods according to the total profit and running time  methods with small-sized set of trip requests.}
\end{figure*}

\noindent\textbf{Greedy-Max-Profit-Based Method}:\\
We now describe a greedy method for profit-maximizing real-time matching. We begin by assigning each of $n$ input users to their ``own'' cab, i.e., no two users are assigned to the same cab.\footnote{We assume that all cabs are initially empty; however, our methods extend to a more general case where some cabs already have passengers in them with their own constraints to begin with.} We maintain two collections of cabs, namely, a set of unavailable cabs $\mathcal{U}$ and a current pool of cabs $\mathcal{A}$. We initialize $\mathcal{A}$ with the $n$ single-user cabs and $\mathcal{U} = \emptyset$. We repeatedly perform the following steps until $\mathcal{A}$ is empty, at which point, we return $\mathcal{U}$ as the output:
\begin{enumerate}[(a)]
\item For every two cabs $j$ and $k$ in $\mathcal{A}$, if it is feasible to \textit{merge} them subject to the capacity constraint (and any other optional constraints), compute the corresponding incremental profit $p_{jk} = \Delta p(\mathcal{S}_j,\mathcal{S}_k)$.
\item Find two cabs $j^*$ and $k^*$ such that $p_{j^*k^*}$ is maximum. If $p_{j^*k^*}<0$, terminate and return $\mathcal{U}\cup \mathcal{A}$ as output. Otherwise, \textit{merge} them and replace the individual cabs $j^*$ and $k^*$ in $\mathcal{A}$ with the merged cab. If the merged cab has $\zeta$ passengers, move it from $\mathcal{A}$ to $\mathcal{U}$.
\end{enumerate}
Note that each iteration either terminates the method, or decreases the size of $\mathcal{A}$ by at least one; therefore, the method terminates after at most $n-1$ iterations. The time complexity of this method is $O(n^2\log n)$.\\

\noindent \textbf{Greedy-Order-Based Methods}:\\
Now we describe a family of faster greedy heuristics (parameterized by an order protocol). We create an ordered list $\mathcal{L}$ of the initial $n$ single-user cabs according to a certain order, and initialize the set of unavailable cabs $\mathcal{U}=\emptyset$. At each iteration, we pick the cab at the top of the list (say $j$) and \textit{merge} it with the first cab $k$ down the list such that $p_{jk}>0$, subject to the capacity constraint (and any other optional constraints), remove the individual cabs $j$ and $k$ from $\mathcal{L}$ and insert the merged cab at the appropriate position in $\mathcal{L}$. If no such cab $k$ exists, then $j$ is moved from $\mathcal{L}$ to $\mathcal{U}$. If the merged cab has $\zeta$ passengers, we move it from $\mathcal{L}$ to $\mathcal{U}$. These methods have a time complexity of $O(n\log n)$.

\subsection{Experimental Performance Evaluation}\label{sec:profit-method-comparison}
In order to compare the performance of our proposed greedy methods, we evaluate them against the optimal solution computed by solving the ILP of Section~\ref{sec:ILP}. Since the ILP takes a prohibitively long time to output the optimal solution, we first evaluate the near-optimality of our methods using a small dataset, and then use a larger dataset for demonstrating the scalability and comparing the speed-optimality trade-off of our proposed heuristics among themselves. We consider two orders when evaluating Greedy-Order-Based methods:
\begin{enumerate}[(a)]
\item Distance Order: Decreasing order of $d(\mathcal{S})$.
\item Profit Order: Increasing order of $p(\mathcal{S})$.\\
\end{enumerate}

\noindent\textbf{Experimental Setup}:\\
We evaluate our methods using publicly available New York City (NYC) taxi trip data~\cite{ny-taxi-data}. The $19$GB dataset has logs for all NYC trips taken in $2013$. Each trip in the dataset specifies pickup and dropoff coordinates, pickup time, dropoff time, trip distance and travel time. Each trip is considered as a ``user'' for our simulation. We assume that each cab serves at most $\zeta=3$ users.

We discretize the underlying NYC map into grids of $100 m^2$ area and select $2543$ representative landmarks using the methodology of~\cite{drs-id}. We precompute the inter-landmark distance and time using Open Trip Planner~\cite{OTP}. For each user, we map their source and destination coordinates to their nearest representative landmarks and use the precomputed values for computing necessary metrics such as profit, detour, driving distance and time.
%
The values $c_b$, $c_d$, $c_t$, and $f_d$ are taken from publicly available Lyft user pricing and driver payment data~\cite{lyftDriver,lyftPrice}.\\

\noindent\textbf{Experimental Results}:\\
We select NYC taxi trips on randomly chosen days of $2013$, between $7$:$45$ - $8$:$00$ pm. We intentionally choose the rush times, since many shared rides would be possible, and it would be interesting to see how our proposed heuristics perform. We run our methods on several subsamples by considering a different number of users each time ($5$, $8$, $10$, $15$, $20$, and $22$). We compare the total profit generated (Figure~\ref{fig:profit_20}) and the total time taken (Figure~\ref{fig:time_20}) across each method, averaged over the small-sized subsamples of users. We observe that the total profit obtained by the Greedy-Max-Profit-Based method is up to $90\%$ optimal and runs $10^{5}$ times faster than the optimal algorithm on a machine with a quad Intel Core i7 processor with $32$GB RAM. The optimal algorithm becomes intractable for instances where the number of users exceeds $25$ (terminated after over $14.22$ hours due to lack of memory).

\begin{figure}
\centering
{\includegraphics[width=0.6\columnwidth]{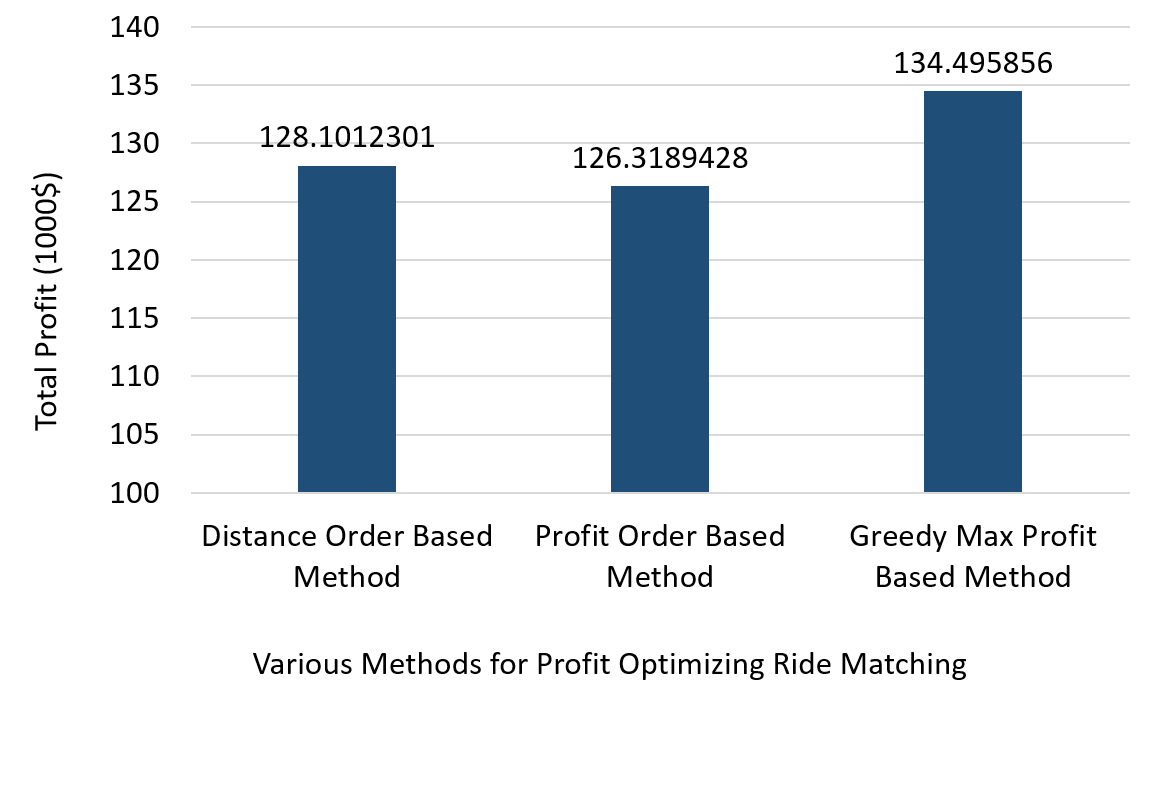}}
\caption{Total Profit earned using various methods on large dataset of {\raise.17ex\hbox{$\scriptstyle\sim$}}$19000$ trip requests.}\label{fig:profit_all}
\end{figure}

Our proposed algorithms scale to a larger input set with {\raise.17ex\hbox{$\scriptstyle\sim$}}$19000$ users (obtained from a one-hour slot). The Greedy-Max-Profit-Based method takes {\raise.17ex\hbox{$\scriptstyle\sim$}}$35$ minutes while the Distance-Order-Based and Profit-Order-Based methods take $7$ and $4$ minutes, respectively.\footnote{If the matching algorithm is invoked once per minute, in NYC, there are only {\raise.17ex\hbox{$\scriptstyle\sim$}}$300$ requests per minute on average, for which all of our proposed greedy methods finish in less than a second.} Figure~\ref{fig:profit_all} shows that the performance of the former is significantly better than that of the latter methods in terms of the total profit. Thus, the speed-optimality trade-off becomes important in choosing the method that best suits the needs of a service provider. Finally, we also observe that the Greedy-Max-Profit-Based method matches {\raise.17ex\hbox{$\scriptstyle\sim$}}$19000$ users into $6701$ cabs (reducing the operational cost significantly), obtaining a profit of {\raise.17ex\hbox{$\scriptstyle\sim$}}$\$134,500$.

\section{Learning the Optimal Discount Policy}\label{sec:learning}

\begin{figure}[ht]
    \centering{\includegraphics[width=0.6\columnwidth]{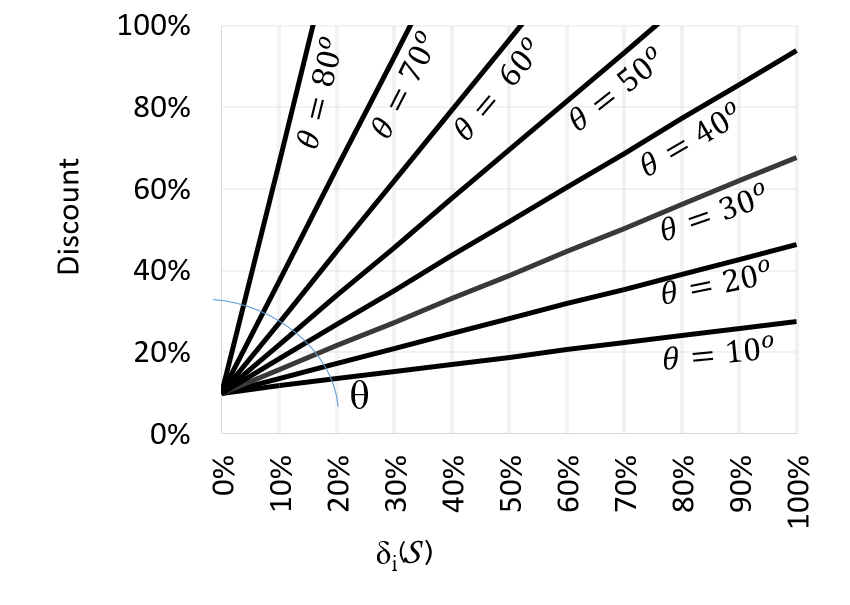}}
\caption{Various discount functions.}
\label{fig:discountFunctions}
\end{figure}

In the previous section, we laid out efficient ride-matching heuristics to maximize profit under a \textit{fixed} detour-based discount for users. In this section, our goal is to learn the optimal detour-based discount policy, under a \textit{fixed} ride-matching algorithm. For simplicity, we assume that the detour-based discount is a linear function of the fractional distance-wise detour. Thus, the discount to a user $i \in \mathcal{S}$ is given by
\begin{equation}\label{eq:discount_function}
f_p^{\theta}(\delta_i(\mathcal{S})) = \tan({\theta}){\delta_i(\mathcal{S})} + b,
\end{equation}
where $b$ is a constant denoting the minimum discount given to a user (to incentivize them to opt for ridesharing in the first place), and the discount parameter $\theta\in[0^{\circ}, 90^{\circ})$ governs how steeply the discount increases with the detour. Figure~\ref{fig:discountFunctions} plots~\eqref{eq:discount_function} for various values of $\theta$, assuming $b=10\%$.

Let $\mathbb{P}(\theta)$ denote the probability of a user opting to rideshare, a nondecreasing function of $\theta$. A larger $\theta$ would lead to a larger population of users who are willing to rideshare, and potentially more shared rides; however, the profit from each shared ride would be smaller. Since the expected total profit depends on the number of users opting to rideshare, as well as the per-ride profit, it is important to find a value of $\theta$ that provides the optimal balance to maximize the expected total profit. We propose a method, which we call \textit{Iterative Discount Function Learning Algorithm} (\textbf{\texttt{IDFLA}}), that learns such a $\theta$ over a period of time.

\begin{algorithm}
\caption{{\textbf{\texttt{IDFLA}} \label{algo:IDFLA}}}
\begin{algorithmic}[1]
{
\State \textbf{Input:} A finite number of $\theta$ values in an array $\Theta$
\State $\mathtt{P}_{h}$ stores the total profit with $\Theta[h]$
\State $k_{h}$ stores the number of days $\Theta[h]$ is declared
\For{$t \gets 1$ \textbf{to} $|\Theta|$}
    \State Choose $\Theta[t]$ on $t^{th}$ day
    \State Observe profit earned $p$;
    \State Set $\mathtt{P}_t\leftarrow p$;
    \State Set $k_t\leftarrow 1$;
\EndFor
\State Find $h^* = \mathop{\arg\max}_h\ \mathtt{P}_h$;
\For{$t \gets |\Theta|+1, |\Theta|+2, \ldots$ }
    \State Choose $\Theta[h^*]$ on $t^{th}$ day
    \State Observe profit earned $p$;
    \State Update $\mathtt{P}_{h^*}\leftarrow \mathtt{P}_{h^*}+p$;
    \State Update $k_{h^*}\leftarrow k_{h^*}+1$;
    \State Find $h^*\!=\! \underset{h}{\mathop{\arg\max}}\ \left(\frac{\mathtt{P}_{h}}{k_{h}}\!\!+\!\!\sqrt{\frac{2\ln{t}}{k_h}}\right)$ for $(t+1)^{th}$ day;
\EndFor
}
\end{algorithmic}
\end{algorithm}

Each day, a value of $\theta$ is declared by the service provider. In response, the users who opt for ridesharing are then matched throughout the day using one of the heuristics from Section~\ref{sec:model}. The total profit earned is computed at the end of the day. \textbf{\texttt{IDFLA}} (Algorithm~\ref{algo:IDFLA}) learns the expected daily profit earned using different values of $\theta$, and eventually converges to the best $\theta^*$ that maximizes this quantity. This kind of sequential decision making in stochastic environment can be categorized as a stochastic multi-armed bandit problem~\cite{bubeck2012regret}, where the $\theta$ values are the arms, and the daily profit earned by using a certain $\theta$ is the reward earned by pulling the corresponding arm. \textbf{\texttt{IDFLA}} uses a popular technique for solving stochastic multi-armed bandit problems, called \textbf{\texttt{UCB1}}~\cite{auer2002finite}.

\subsection{Performance Evaluation of IDFLA}

\begin{figure}
\centering
{\includegraphics[width=0.6\columnwidth]{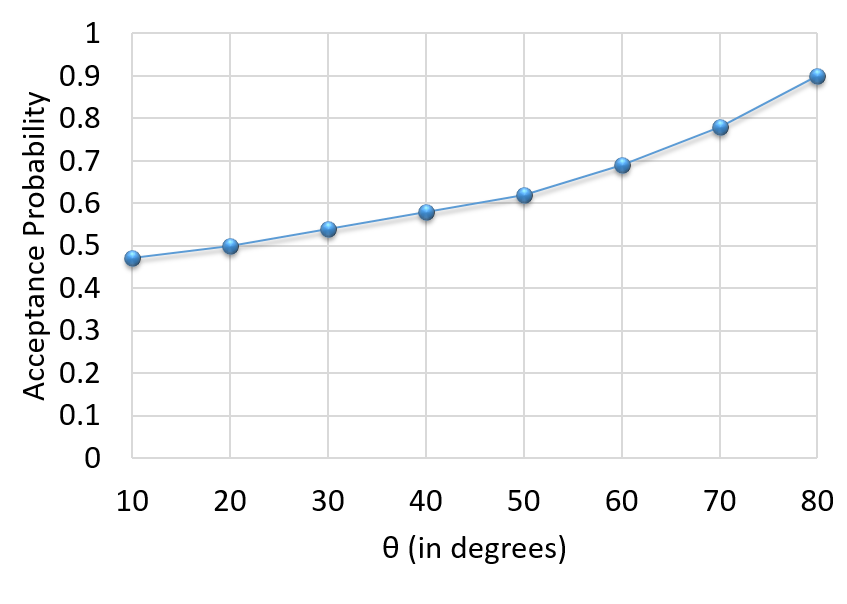}}
\caption{Probability of opting for a shared ride.}
\label{fig:acceptance_probability}
\end{figure}

We evaluate our learning method \textbf{\texttt{IDFLA}} under the Greedy-Max-Profit-Based ride-matching heuristic. We take $\Theta = \{10^{\circ}, 20^{\circ}, 30^{\circ}, 40^{\circ}, 50^{\circ}, 60^{\circ}, 70^{\circ}, 80^{\circ}\}$, and run \textbf{\texttt{IDFLA}} to find near-optimal discount parameter, $\Theta[h^*]$. The value of optimal  discount parameter depends on $\mathbb{P}(\cdot)$ which is the probability of opting rideshare. For simulating the user-behavior, we need to assume $\mathbb{P}(\cdot)$ to be any non-decreasing function of $\theta$. In practice, this probability $\mathbb{P}(\theta)$ is estimated while the optimal $\theta$ is learnt (step 14 of \textbf{\texttt{IDFLA}}). This corresponds to the fraction of times ($k_{\theta}$) the users have opted for ridesharing whenever the discount parameter was $\theta$. The $\mathbb{P}(\theta)$ function that we have assumed for our simulation, is shown in Figure~\ref{fig:acceptance_probability}.

\begin{figure}
    \centering
    \includegraphics[width=0.6\columnwidth]{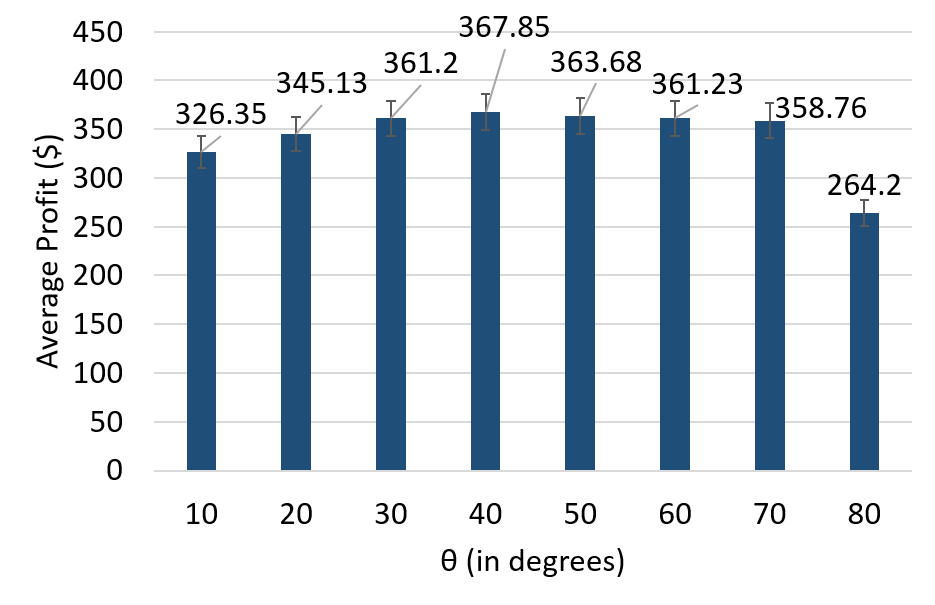}
    \caption{The average daily profit obtained by various discount functions over 52 days}
    \label{fig:bdf_52}
\end{figure}

Taxi trips in NYC between $7$:$45$ - $8$:$00$ pm of the $52$ Wednesdays of $2013$
were used to compare the performance of \textbf{\texttt{IDFLA}} with an ``oracle'' method (\textbf{\texttt{BDF}}) that knows the Best Discount Function apriori. Figure~\ref{fig:bdf_52} shows the average daily profit obtained by declaring the same discount function everyday. We observe that the best discount function is obtained when $\theta = 40^{\circ}$; therefore, \textbf{\texttt{BDF}} chooses $\theta = 40^{\circ}$ everyday.

\begin{figure}
    \centering
    \subfloat[][{Comparing the daily profit obtained by \textbf{\texttt{IDFLA}} and \textbf{\texttt{BDF}}.}\label{fig:total_profit_52}]{\includegraphics[width=0.6\columnwidth]{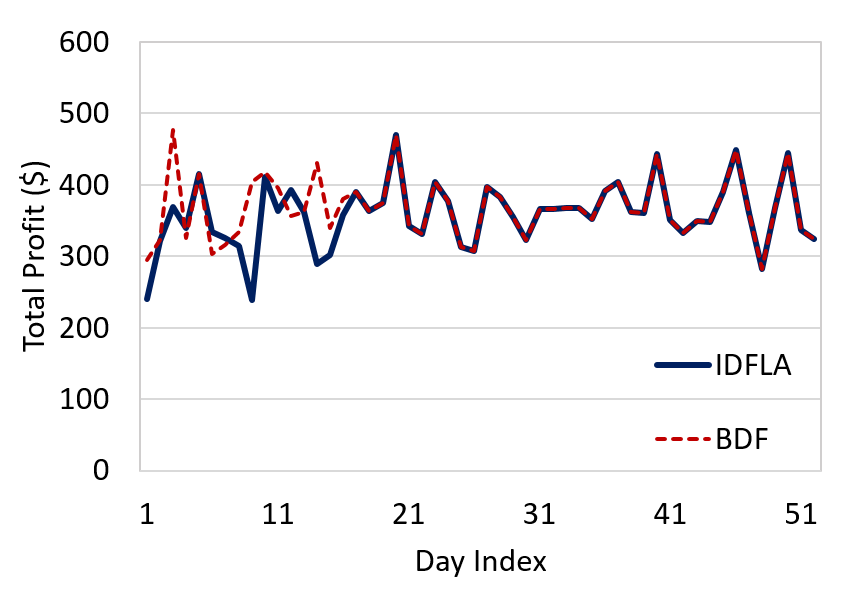}}\\
    \subfloat[][{The average daily profit of \textbf{\texttt{IDFLA}} converging to that of \textbf{\texttt{BDF}}.}\label{fig:average_profit_52}]{\includegraphics[width=0.6\columnwidth]{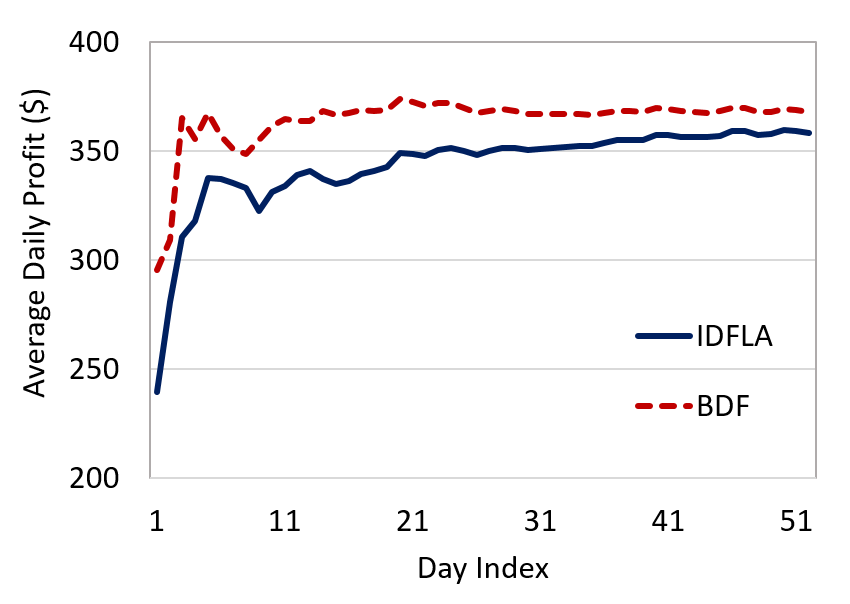}}
    \caption{Comparing \textbf{\texttt{IDFLA}} and \textbf{\texttt{BDF}}}\label{fig:52}
\end{figure}

Figure~\ref{fig:total_profit_52} shows that the discount parameter selected by \textbf{\texttt{IDFLA}} converges to that of \textbf{\texttt{BDF}} in $16$ days. We also observe (Figure~\ref{fig:average_profit_52}) that the average daily profit, that is, the cumulative total profit over $t$ days divided by $t$, converges to that of the best discount function \textbf{\texttt{BDF}}. We show that after $52$ days, the difference\footnote{The difference between the average reward of a learning algorithm to that of the optimal is called $\textit{regret}$. An upper bound on the regret of \textbf{\texttt{UCB1}} after $T$ iterations is $O(\log T)$~\cite{auer2002finite}, when each arm's reward is a bounded random variable.} in average daily profit obtained by \texttt{\textbf{IDFLA}} and \texttt{\textbf{BDF}} is just $\$9.58$, or $2.67\%$. Thus, \textbf{\texttt{IDFLA}} quickly learns the optimal discount parameter that can help the service provider maximize profit while keeping users sufficiently happy by providing suitable discounts.

\section{Impact of Detour-Aware Routing}\label{sec:qos-impact}

In conjunction with desirable pricing schemes for ridesharing users, recent work has proposed detour-aware routing policies that impose upper bounds on the total detour~\cite{Jung2013} (static) or incremental detours~\cite{DBLP:journals/corr/GopalakrishnanM16} (dynamic). Such detour-aware routing policies work together with detour-based discount policies in order to enhance the quality of the ride experience for ridesharing users. This induces a similar trade-off as the one in Section~\ref{sec:learning}: imposing additional quality-enhancing constraints into the ride-matching process incentivizes greater adoption of ridesharing and increases the market share of the ridesharing population, but could result in reduced average per-ride profit. In this section, we consider the impact of one such detour-aware routing policy, based on the concept of \textit{Sequential Individual Rationality} (SIR), introduced by~\cite{DBLP:journals/corr/GopalakrishnanM16}.

SIR guarantees that the \textit{disutility} to existing users (sum of the monetary cost and an ``inconvenience cost'' due to detours) in a shared ride is non-increasing as additional users are picked up. A parameterized version of SIR, called SIR-$\gamma$, ensures that the incremental benefit (decrease in disutility) to an existing user upon picking up a new user is at least $\gamma$, where $\gamma\geq 0$. Imposing SIR-$\gamma$ results in a reduced set of feasible matched rides, leading not only to reduced average per-ride profit, but also a potential loss due to the ``up-front discounts'' (the parameter $b$ in~\eqref{eq:discount_function}) given to the passengers who opt for ridesharing, but end up being unmatched. However, the resulting detour-aware routing policy improves the quality of shared rides, thereby increasing the market share of ridesharing users, which could offset these negative effects. Unlike~\cite{DBLP:journals/corr/GopalakrishnanM16} that looks at the detour-aware routing policy from the user perspective (especially fairness), in this work we answer such a scheme's impact on the service provider's profit.

We investigate this phenomenon by performing experiments to address the following questions:
\begin{enumerate}[(a)]
\item What is the trade-off between fewer matched rides and increased market share, compared to a scenario where no SIR-$\gamma$ is imposed? How does this affect the profit?
\item How sensitive is the profit to the fraction of additional users who opt for ridesharing in response to the adoption of detour-aware routing by imposing SIR-$\gamma$?
\end{enumerate}

\subsection{Experimental Setup}
As before, we use NYC taxi trip data for our experiments. We randomly select a weekday of $2013$ and a one hour time slot on that day. We then split the hour into $60$ one-minute instances, and collect the sources and destinations for the users who initiated trips in each instance. The average number of users per instance is $210.85$. While reporting our results, we take an average across these $60$ instances.

The parameters $c_b$, $c_d$, $c_t$, and $f_d$ are taken from publicly available Lyft user pricing and driver payment for NYC~\cite{lyftDriver,lyftPrice}. We select ``detour-sensitivities'' $\alpha_i$ (higher values imply more aversion to detours) randomly from $[\$0,\$5]$ per mile.
We assume that each cab serves at most $\zeta=2$ users, so that the profit maximization problem can be solved optimally in polynomial time using Edmond's algorithm~\cite{edmonds1965paths}.\footnote{The experiments can be carried out for $\zeta>2$ as well, by using one of our proposed heuristics from Section~\ref{sec:model} instead.}
The driving distance and time between locations are obtained by querying Open Trip Planner~\cite{OTP}. 
The parameters of the linear detour-based discount policy~\eqref{eq:discount_function} are $\theta=40^{\circ}$ and $b=10\%$. To begin with, before imposing SIR-$\gamma$, we assume that the service provider's initial market share is $60\%$, within which each user opts for ridesharing with a probability directly proportional to $\theta$ and $b$, and inversely proportional to $\alpha_i$.

For each instance, we keep the following fixed: user origins, destinations, whether they are part of the initial market share, and if so, whether they opted for ridesharing. Then, we generate $30$ realizations, where each user reacts to the adoption of a detour-aware routing policy (SIR-$\gamma$), by flipping biased coins parameterized by $c_{\textrm{in}}$, which is a measure of how strongly users value a quality-enhancing detour-aware routing policy. This participatory behavior is one of many that can be used. For any behavior model, our primary interest is the fraction of users who end up opting for ridesharing, and not the generative process itself. In the coin-flip-based process, the first group of users, who are in the initial market share but did not opt for ridesharing, flip with bias $p_{\textrm{in}}$ that is directly proportional to $\theta$, $b$, and $\gamma$, and inversely proportional to $\alpha_i$, with $c_{\textrm{in}}$ denoting the constant of proportionality. The second group of users, who are outside the initial market share, flip with bias $0.5p_{\textrm{in}}$ (the factor $0.5$ is hand-picked). The users who had already opted for ridesharing remain so. For consistency, we ensure that as we increase the value of $\gamma$, only users who have previously not opted for ridesharing flip again.

We calculate the profit obtained by matching users in the new market share. The parameter $c_{\textrm{in}}$ is varied from $100$ and $1000$, while $\gamma$ is varied from $0$ to $0.9$.

\subsection{Experimental Results}

\begin{figure}
\centering
    \includegraphics[width=0.6\columnwidth]{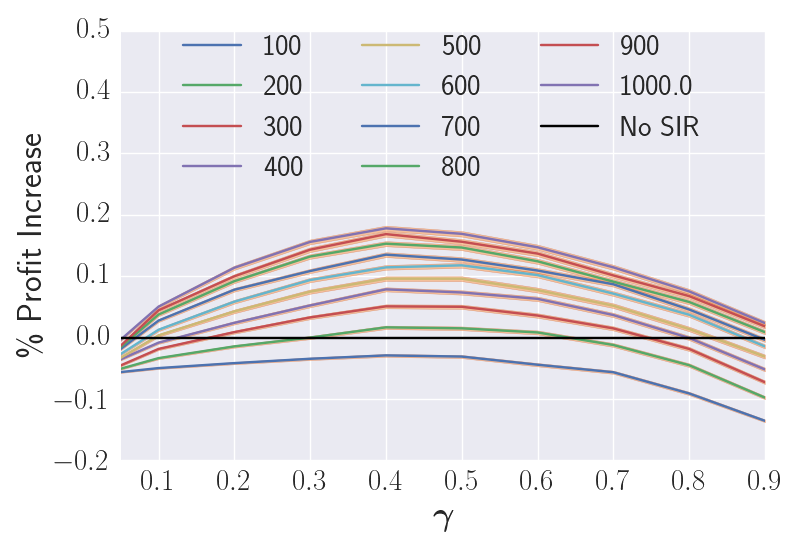}
    \caption{Profit versus SIR parameter $\gamma$.\label{fig:profit-gamma}}
\end{figure}

Figure~\ref{fig:profit-gamma} shows profitability as $\gamma$ increases for various values of $c_{\textrm{in}}$. The baseline profit, without SIR-$\gamma$, is shown as a horizontal line. For any value of $c_{\textrm{in}}\geq 200$, there is a range of $\gamma$ for which a higher profit can be obtained by imposing the appropriate SIR-$\gamma$ policy. Moreover, when $\gamma$ is high, the profit declines, because too many users remain unmatched.

\begin{figure}
\centering
    \includegraphics[width=0.6\columnwidth]{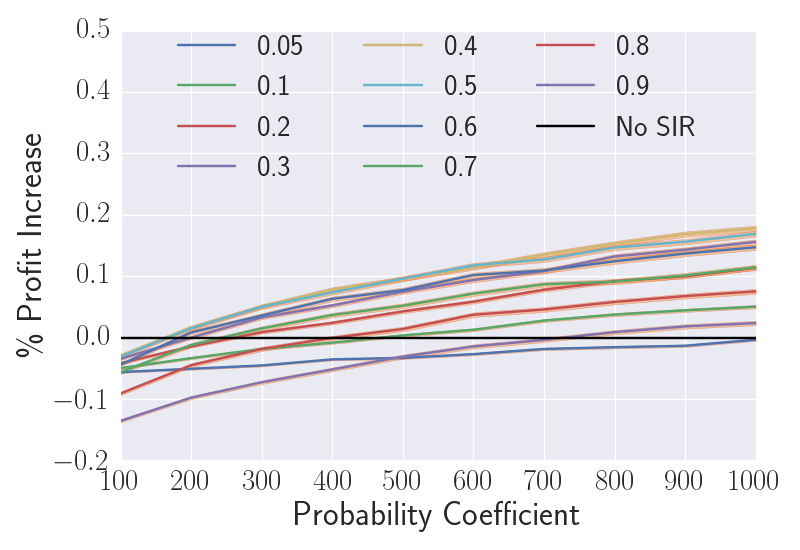}
    \caption{Sensitivity of profit to $c_{\textrm{in}}$.\label{fig:profit-probability}}
\end{figure}

Figure~\ref{fig:profit-probability} shows a different visualization by plotting profitability as a function of $c_{\textrm{in}}$, for different values of $\gamma$.

\begin{figure}
\centering
  \includegraphics[width=0.6\columnwidth]{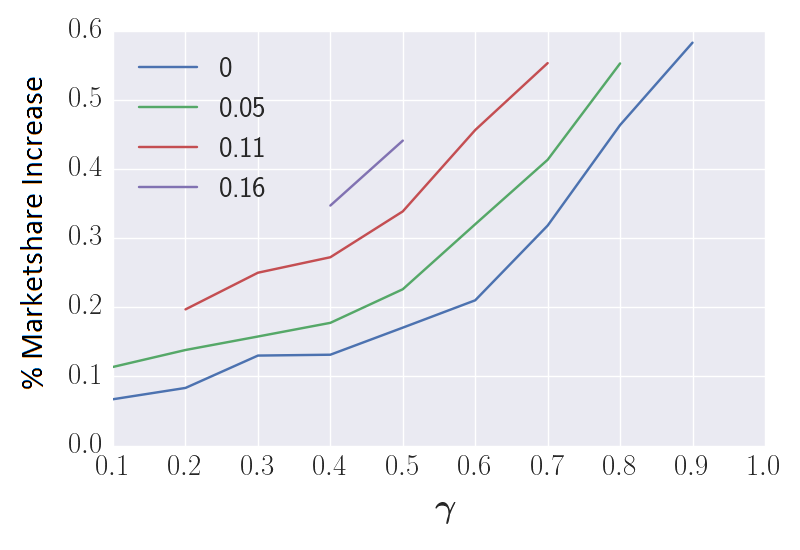}
    \caption{Marketshare versus $\gamma$ for a given profit lower bound.\label{fig:ms-gamma}}
\end{figure}

\begin{figure}
\centering
    \includegraphics[width=0.6\columnwidth]{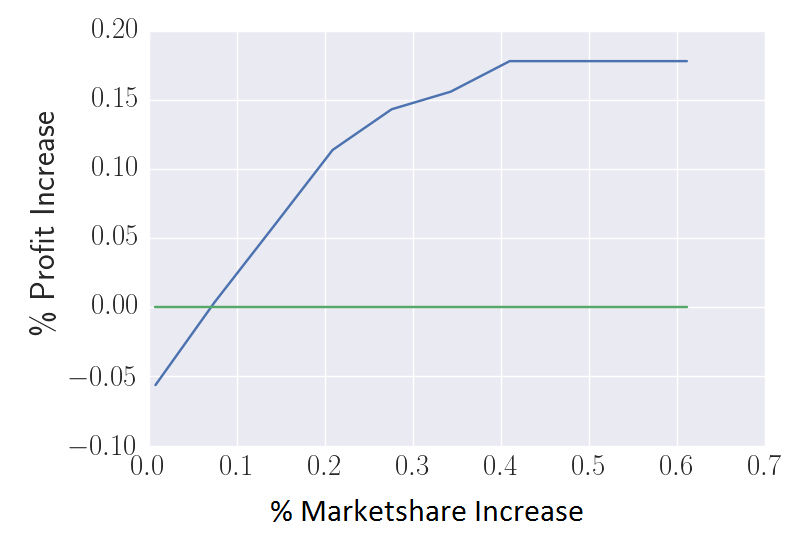}
    \caption{Profit versus marketshare increase.\label{fig:profit-ms}}
\end{figure}

Figure~\ref{fig:ms-gamma} plots the minimum increase in market share that is necessary to meet a given lower bound on the profit, as a function of $\gamma$. 
For instance, to be $5\%$ more profitable, a $13\%$ increase is needed, as shown by the green curve. Alternatively, Figure~\ref{fig:profit-ms} shows the maximum possible increase in profit (over all $\gamma$), as a function of the increase in market share. The key takeaway here is that a $7\%$ increase\footnote{A $7\%$ increase in market share in response to a quality-enhancing policy is not unrealistic. Lyft's market share increased by over $70\%$ in $2016$, owing in large part to better quality of service, resulting in Lyft's customers being happier than those of their primary competitor, Uber~\cite{lyftMarket}.} in market share is sufficient to recover from the negative effects (on the profit) of adopting a detour-aware routing policy.

\section{Concluding Remarks}

In this paper, we adopt a profit-centric view of real-time ridesharing system design, and undertake an empirical investigation of three key elements that impact profit, namely (a)~real-time ride matching, (b)~user discounts, and (c)~detour-aware routing, by using publicly available NYC taxi trip data.

An optimal ride matching algorithm, even on a small set of $22$ requests, can take hours, while our greedy heuristics generate near-optimal profit in less than $1/10^5$ of the time, and scale well, e.g., matching $18794$ trip requests into $6701$ cabs, while offering different trade-offs between optimality and running time. Our adaptive learning algorithm called \textbf{\texttt{IDFLA}} learns the optimal detour-based discount parameter that achieves the right balance between the number of passengers who opt for ridesharing and the portion of savings taken as profit from each shared ride, within $16$ iterations. Finally, we demonstrate that even a small market share increase of $7\%$, in response to other restricted detour-aware policies such as sequential individual rationality, is sufficient to counter the reduced per-ride profit due to the restricted constraints. Thus, one of our contributions is also a holistic treatment of the framework for detour-aware policies.

Given the encouraging empirical results of our proposed solutions for profit optimization, future work should explore corresponding theoretical results, e.g., approximation guarantees for the proposed greedy ride-matching heuristics, regret bounds for the proposed learning algorithm \textbf{\texttt{IDFLA}}, and a formal mathematical analysis to characterize the dependence between profit under a detour-aware routing policy and market share of users who opt for ridesharing. Another interesting direction to explore is the interplay between these different elements, e.g., how does a detour-aware routing policy affect the performance of the ride-matching heuristics and the convergence of the learning algorithm?

\bibliography{profitOptimization}
\bibliographystyle{myIEEEtran}

\end{document}